\input amstex
\loadbold
\documentstyle{amsppt}
\magnification=1200

\pagewidth{5.75truein}
\pageheight{9.2truein}
\hcorrection{0.8truecm}

\

\vskip -5mm

%%{\eightpoint \noindent {\it Open
%%Problems. Workshop on Group Actions on
%%Rational Varieties,
%%
%%\noindent McGill
%%Univ.~and~Univ.~of~Montreal, Canada,
%%March 2002}.
%%
%%\noindent G.\;Freudenburg, P.\;Russell
%%Eds., Reports Centre Interuniv.
%%
%%\noindent
%% Calcul. Math.
%%Alg\'ebrique. Concordia, Laval, McGill,
%%2002--04, 13--14. }
%%
%%\vskip 2mm

%%{\eightpoint \noindent To appear in:}

%%\vskip 2mm

{\eightpoint \noindent {\it Open Problems in Affine
Algebraic Geometry},

\noindent in: {\it Affine Algebraic Geometry},

\noindent J.~Gutierrez, V.~Shpilrain, J.--T.~Yu,
Eds.,

\noindent Contemp. Math., Amer. Math. Soc.

\noindent To appear.\footnote{Preprint: Open
problems session of {\it Workshop on Group Actions
on Rational Varieties}, McGill Univ. and Univ. of
Montreal, Canada, March 2002 (G.~Freudenburg,
P.~Russell, Eds.), Reports Centre Interuniv.,
Calcul. Math. Alg\'ebr., Concordia, Laval, McGill,
2002--04, 13--14.}

\vskip 2cm

\topmatter
\title
Problems for problem session
\endtitle
\author
Vladimir L. Popov
\endauthor
\address
Steklov Mathematical Institute, Russian
Academy of Sciences, Gubkina 8, Moscow,
117966, Russia
\endaddress
\email popovvl\@orc.ru
\endemail

%%\date
%%March 20, 2002
%%\enddate
\endtopmatter
\document

Below are the problems that I formulated at Open
Problems Session of {\it Workshop on Group Actions
on Rational Varieties}, McGill University and
University of Montreal, Canada, March 2002.

\medskip

Let $k$ be an algebraically closed field
of characteristic $0$.

\medskip

 {\bf 1. \ Roots of the Affine Cremona
Group}

\medskip

Let $k^{[n]}:=k[x_1,\ldots,x_n]$ be the
polynomial algebra in variables
$x_1,\ldots,x_n$ over $k$, and let $n>1$.
Put
$$\textstyle D_i:=\frac{\partial}{\partial
x_i}.$$ {\it The affine Cremona group}
$\operatorname{Aut}_k k^{[n]}$ is an
infinite dimensional algebraic group,
 and the set of
`volume preserving' transformations,
$$\textstyle \operatorname{Aut}^\ast_k
k^{[n]}:= \{\sigma\in \operatorname{Aut}_k
k^{[n]}\mid \det\bigl(
%%\frac{\partial
D_i(\sigma(x_j))
%%}{\partial x_j}
\bigr)_{1\leq i, j\leq n}=1\}, $$ is its
closed normal subgroup, \cite{Sh},
\cite{Ka}. The group
$\operatorname{Aut}^\ast_k k^{[n]} $ is an
infinite dimensional simple algebraic
group, \cite{Sh}. If $D$ is a locally
nilpotent $k$-derivation of $k^{[n]}$,
then it is easily seen that $\exp tD\in
\operatorname{Aut}^\ast_k k^{[n]}$ for any
$t\in k$, so $D$ lies in
$\operatorname{Lie}
(\operatorname{Aut}^\ast_k k^{[n]} )$.

It follows from \cite{BB1}, \cite{BB2}
that $n-1$ is the ma\-xi\-mum of
dimensions of algebraic tori contained in
$\operatorname{Aut}^\ast_k k^{[n]}$, and
that every algebraic torus in
$\operatorname{Aut}^\ast_k k^{[n]} $ of
dimension $n-1$ is conjugate to the
`diagonal' torus
$$\textstyle T:=\{\sigma:=
\operatorname{diag}(t_1,\ldots,t_n)\in
\operatorname{Aut}^\ast_k k^{[n]} \mid
\sigma (x_i)=t_ix_i, \ t_i\in k,\ \prod
_{i=1}^n t_i=1\}.$$

 \vskip 1.5mm

Mimicing the notion of root
%%from
%%classical
from the theory of finite dimensional
linear algebraic groups, cf.,
e.g.,\;\cite{Sp}, introduce the following

\proclaim{Definition 1.1} {\rm A nonzero
locally nilpotent derivation $D$ of
$k^{[n]}$ is called a} root vector {\rm of
$\operatorname{Aut}^\ast_k k^{[n]} $ with
respect to $T$ if there is a nontrivial
character $\chi \in \operatorname{Hom}(T,
{\bold G}_m)$ such that}
$$\sigma\circ
D\circ\sigma^{-1}=\chi(\sigma) D \ \
\forall \ \sigma\in T.$$

\noindent {\rm The character $\chi$ is
called the} root {\rm of
$\operatorname{Aut}^\ast_kk^{[n]}$ with
respect to $T$ corresponding to $D$}.
\endproclaim

\proclaim{Problem 1.1}{\it Find all roots
and root vectors of
$\operatorname{Aut}^\ast_kk^{[n]}$ with
respect to $T$.}
\endproclaim

Among the root vectors there are evident
`classical' ones
$$
\textstyle x_iD_j,\ \text{where}\ i\neq j,
\ \text{and} \ D_i
$$
contained respectively in the Lie algebras
of the special linear group
$\operatorname{SL}_n$ (acting naturally on
$k^{[n]}$ by linear transformations of
$x_1,\ldots, x_n$) and the group of
translations
$$\operatorname{T}_n:=\{(x_1,\ldots,x_n)\mapsto
(x_1+c_1,\ldots, x_n+c_n)\mid c_i\in
k\}.$$
 The
correspon\-ding roots are respectively
$\varepsilon_i\varepsilon^{-1}_j$ and
$\varepsilon_i$ where
$$
\varepsilon_i(
\operatorname{diag}(t_1,\ldots, t_n))=t_i.
$$

\proclaim{Problem 1.2}{\it Describe the
centralizer and the normalizer of $\,
T$\,in
$\operatorname{Aut}^\ast_kk^{[n]}$.}
\endproclaim

\noindent{\it Remark $1.1$.} The question
whether $\operatorname{Aut}^\ast_kk^{[n]}$
has no maximal (with respect to inclusion)
algebraic tori of dimension $< n-1$ is
equivalent to asking whether every `volume
preserving' algebraic torus action on
$\bold A^n$ is linearizable. The first
nontrivial case is $n=3$. Putting together
the results of \cite{KR}, \cite{KR1},
\cite{Ko}, \cite{KM}, \cite{KR2},
\cite{KR3} (see a historical account in
\cite{KMKR}) one obtains that in this case
the answer is affirmative.

\bigskip

{\bf 2. \ Generators of
$\operatorname{Aut}^\ast_kk^{[n]}$}

\medskip

\proclaim{Definition 2.1} {\rm A subgroup
$G$ of $\operatorname{Aut}^*_kk^{[n]}$ is
called $\partial$-{\it generated}  if
there exists a set $\Cal G$ of locally
nilpotent derivations of $k^{[n]}$ such
that $G$ is generated by the elements
$\exp D$, where $D\in\Cal G$.}
\endproclaim

Recall (cf.~\cite{P}) that, given a
locally nilpotent derivation $D$ of
$k^{[n]}$, there is a simple method of
constructing new such derivations: if $f$
is an element of
$\operatorname{Ker}D:=\{h\in k^{[n]}\mid
Dh=0\}$, then $f\hskip -.2mm D$ is a
locally nilpotent derivation of $k^{[n]}$
as well. Thus for any such $f$ we obtain
an element $\exp f\hskip -.2mmD$ of
$\operatorname{Aut}^*_kk^{[n]}$.

\proclaim{Definition 2.1} {\rm A subgroup
$G$ of $\operatorname{Aut}^*_kk^{[n]}$ is
called {\it finitely $\partial$-generated}
if there exists a finite set $\Cal G$ of
locally nilpotent derivations of $k^{[n]}$
such that $G$ is generated by the elements
$\exp f\hskip -.2mmD$ where $D\in\Cal G$
and $f\in \operatorname{Ker} D$.}
\endproclaim

\noindent {\it Example $2.1$.} Since the
special linear group $\operatorname{SL}_n$
is generated by the root subgroups $\{\exp
tx_i D_j \mid t\in k\}$, $1\leq i\neq
j\leq n$, cf.\,\cite{Sp}, it is  finitely
$\partial$-generated (with $\Cal G =\{D_1,
\ldots, D_n \}$). By the same reason, any
connected semisimple algebraic subgroup
$G$ of $\operatorname{Aut}^*_kk^{[n]}$ is
finitely $\partial$-generated. On the
other hand, any algebraic torus in
$\operatorname{Aut}^*_kk^{[n]}$ is not
$\partial$-generated.

\medskip

\noindent {\it Example $2.2$.} Since the
automorphism
$$(x_1,\ldots, x_n)\mapsto
(x_1+c_1,\ldots, x_n+c_n),\ c_i\in k,$$
coincides with $(\exp
c_1D_1)\circ\ldots\circ (\exp c_nD_n)$,
the group of translations
$\operatorname{T}_n$ is finitely
$\partial$-generated (with $\Cal
G=\{D_1,\ldots, D_n\}$).

\medskip

\noindent{\it Example $2.3$.} Let $n>1$.
If $f_i\in k[x_{i+1},\ldots, x_n]$,
$i=1,\ldots, n-1$, then $(\exp f_{1}D_{1})
\circ\ldots\circ(\exp f_{n-1}D_{n-1})$ is
the automorphism
$$(x_1,\ldots,x_n)\mapsto (x_1+f_1,
x_2+f_2, \ldots, x_{n-1}+f_{n-1}, x_n).
$$
Hence the de Jonqui\`ere group
$\operatorname{J}_n$ consis\-ting of all
such automorphisms is finitely
$\partial$-generated (with $\Cal G
=\{D_1,\ldots, D_{n-1}\}$).

\proclaim{Problem 2.1} Is the group
$\operatorname{Aut}^*_kk^{[n]}\hskip .5mm$
$\partial$-generated? If yes, is it
finitely $\partial$-generated?
\endproclaim

\noindent{\it Example $2.4$}. Since
$\operatorname{Aut}^*_kk^{[1]}=
\operatorname{T}_1$, it follows from
Example 2.1 that the group
$\operatorname{Aut}^*_kk^{[1]}$ is
finitely $\partial$-generated (with $\Cal
G=\{D_1\}$).

\medskip

\noindent{\it Example $2.5$}. Since
$\operatorname{Aut}^*_k k^{[2]}$ is
generated by $ \operatorname{SL}_2$,
$\operatorname{T}_2$ and
$\operatorname{J}_2$, cf.\,\cite{vdK}, it
follows from Examples 2.1, 2.2 and 2.3
that the group
$\operatorname{Aut}^*_kk^{[2]}$ is
finitely $\partial$-generated (with $\Cal
G=\{D_1,D_2\}$).

\medskip

\noindent{\it Remark $2.1$.} Let
$\sigma\in \operatorname{Aut}^*_k k^{[3]}$
be the Nagata automorphism, \cite{N},
$$
(x_1, x_2, x_3)\mapsto
(x_1-2x_2(x_3x_1+x_2^2)-x_3(x_3x_1+x_2^2)^2,
x_2+x_3(x_3x_1+x_2^2), x_3).
$$
By \cite{SU1}, \cite{SU2}, it is not tame,
i.e., does not lie in the subgroup of
$\operatorname{Aut}^*_k k^{[3]}$
ge\-ne\-rated by $\operatorname{SL}_3$,
$\operatorname{T}_3$ and
$\operatorname{J}_3$.  However notice that
$\sigma=\exp f\hskip -.2mmD$, where $D$ is
the locally nilpotent derivation of
$k^{[3]}$ given~by $ D:=-2x_2D_1+x_3D_2$,
and $f:=x_3x_1+x_2^2\in
\operatorname{Ker}D$.

\bigskip

{\bf 3. \ Two Locally Nilpotent
Derivations}

\medskip

Let $P$ and $Q$ be locally nilpotent
derivations of $k^{[n]}$.

\proclaim{Problem 3.1} When is the minimal
closed subgroup of
$\operatorname{Aut}_kk^{[n]}$ containing
the one-dimensional subgroups
$$\{\exp\,tP\mid t\in k\}\hskip 1.5mm
\text{and}\hskip 1.5mm \{\exp\,tQ\mid t\in
k\}
$$ finite dimensional?
\endproclaim

\bigskip

{\bf 4. \ Rationality of Some Homogeneous
Spaces}

\medskip

Let $G$ be a connected semisimple
algebraic group over $k$, and let $\sigma
\in\operatorname{Aut}_k G$ be an element
of finite order. Consider the subgroup
$$G^{\sigma}:=\{g\in G\mid \sigma(g)=g\}.$$

\proclaim{Conjecture 4.1}{\it The
algebraic variety $G/G^{\sigma}$ is
rational.}
\endproclaim

If $\operatorname{ord}\sigma=2$, then, by
\cite{V}, the homogeneous space
$G/G^{\sigma}$ contains a dense open
orbit of a Borel subgroup of $G$. Since
any orbit of a connected sol\-vable
linear algebraic transformation group is
rational, this shows that the conjecture
is true for $\operatorname{ord}\sigma=2$.

\enddocument